\input amstex
\documentstyle{amsppt}
\NoRunningHeads

\magnification\magstep1
\def\s{\vskip6pt}

\def\mpr#1{\;\smash{\mathop{\hbox to 20pt{\rightarrowfill}}\limits^{#1}}\;}
\def\qed{$\;\sqcap\hskip-6.5pt\sqcup$\hfill}

\def\x{{^{_\bullet}}}
\def\Pt{{\Bbb P}}
\def\Ct{{\Bbb C}}
\def\Qt{{\Bbb Q}}
\def\Ft{{\Bbb F}}

\topmatter
\title Formality of equivariant intersection cohomology
of algebraic varieties\endtitle
\author Andrzej Weber \endauthor
\thanks Partially supported by KBN 2P03A 00218 grant and by the
European Commission RTN {\it Geometric Analysis}\endthanks

\affil Institute of Mathematics, Warsaw University \endaffil
\address ul. Banacha 2, 02--097 Warszawa, Poland \endaddress
\email aweber\@mimuw.edu.pl \endemail

\abstract
We present a proof that the equivariant intersection cohomology
of any complete algebraic variety acted by a connected algebraic
group $G$ is a free module over $H^*(BG)$.\endabstract

\endtopmatter
\document

\head Introduction\endhead

Let $X$ be a complex algebraic variety and $G$ be a
linear algebraic group acting algebraically on $X$. Denote by
$IH_G^*(X)=H^*(EG\times_GX;\Cal{ IC}_G^\x)$ the
equivariant intersection cohomology of $X$ with
rational coefficients, first introduced in \cite{Br}. It is a module
over $H^*(BG)$. We say that
$X$ is formal if $IH_G^*(X)$ is free.
We have a spectral sequence
$$E^{pq}_2=H^p(BG,IH^q(X))\Rightarrow IH^{p+q}_G(X)\,.\tag *$$
Suppose that $G$ is connected. Then $BG$ is 1-connected and the
coefficients are not twisted.
Then the spectral sequence (*) degenerates if and only if $X$ is formal.

The aim of this note is to show that
complete $G$-varieties are formal.
We will prove:
\s
\proclaim{ Theorem 1} Let $X$ be a complete algebraic variety acted
by a connected algebraic linear group. Then $X$ is formal.\endproclaim
In this way we add to the list of \cite{GKM, 14.1} an important
class of examples. In fact Theorem 1 is known for specialists. Each
of the proofs comes back to the fundamental paper of
Beilinson-Bernstein-Deligne-Gabber, \cite{BBD}.
After all, in the most general case, Theorem 1 follows from a simple
weight argument, as presented in Remark 7. Nevertheless, we think it is
worth to compare several approaches to equivariant intersection
cohomology and ways of possible proofs. 

I would like to thank Prof.\ Andrzej Bia\l ynicki-Birula for
giving me a hint in a crucial point. I thank Prof.\ Pierre
Deligne for valuable remarks and for assuring me that the methods of
\cite{BBD} work
for algebraic spaces. I thank Prof.\ Frances Kirwan for her comments. I would
like also to thank Marcin Cha\l upnik for help during the preparation of this
note.

\head 1. Proof for projective varieties\endhead

Our proof of Theorem 1
uses the decomposition theorem or reduction to a finite
characteristic, which are strong tools. Thus we think that it is
worth to note that for projective $X$ the situation is much
simpler. We proceed as follows:
\s
\noindent Step 1: We choose an equivariant embedding of $X$ into a
projective space. This is possible by a result of Sumihiro \cite{Sum,
Theorem 1}.
\s
\noindent Step 2: Let $h\in H^2(X)$ be the class of the hyperplane
section. We show that $h$ can be lifted to $H^2_G(X)$. This is so
for $X=\Pt^N$ and $G=GL(N+1)$ acting naturally. Indeed, take a
infinite Grassmannian as a model of $BG$, then
$$EG\times_G\Pt^N=\{(V,l)\,:\, l\subset V\subset
\Ct^\infty\,,\;\dim l=1\,,\;\dim V=N+1\,\}\,.$$ One has a
tautological bundle, its Chern class is a lift of $-h$. If
$X\subset \Pt^N$ and $G$ acts linearly, we construct a lift for
$\Pt^N$ by naturality with respect to $G$ and we restrict it to
$X$.
\s
\noindent Step 3: We have a map of spectral sequences
$$L:E^{p,q}_r\rightarrow E^{p,q+2}_r\,,$$ which is induced by the
multiplication by a cycle (or a form) representing $h$.
Let $d=\dim X$. By the hard Lefschetz theorem
 $$L^k:E^{p,d-k}_2=H^p(BG)\otimes IH^{d-k}(X)\rightarrow
E^{p,d+k}_2=H^p(BG)\otimes IH^{d+k}(X)$$ is an isomorphism for
every $p$ and $k$. From the Lefschetz criterion, \cite{D1}, it follows
that the spectral sequence degenerates, as noticed \cite{Br},
4.2.2.
\s
This argument works for varieties over any field of any
characteristic. One should replace the ordinary cohomology and
intersection cohomology by their \'etale version.
In the Step 3
an isomorphism $IH_G^*(X)\simeq H^*(BG)\otimes IH^*(X)$ can be made
canonical by \cite{D3}. It would depend only on the embedding
$X\subset\Pt^N$ and the representation $G\rightarrow GL(k)$.
\s

For smooth $X$ the Step 2 can also be carried out by methods of
symplectic geometry. Let 
$K$ be a maximal compact subgroup of $G$. Then
$H^*_G(\Pt^N)\simeq H^*_K(\Pt^N)$. The class $h$ is represented by a
$K$-equivariant symplectic form $\omega$. One can define a moment
map by an explicit formula \cite{K, 2.5}. The choice of the moment map
is precisely a lift of $\omega$ to $H^2_K(\Pt^N)$ in the Cartan model
of $\Omega^*(EK\times_K\Pt^N)$, see \cite{AB, 6.18}.
Another proof for smooth projective varieties, can be
obtained by means of a Morse function coming from a moment map, see
\cite{K, 5.8}.
\s
In the rest of the paper we will show how to deduce formality for
an arbitrary complete variety, not necessarily projective.

\head 2. Algebraic models of $EG\times_G X$ \endhead

There are several constructions of $BG$ that can be performed in the
category of algebraic varieties. Of course $BG$ itself cannot be an
algebraic variety, since it cannot be of a finite dimension. 
There are two well known approaches of giving a meaning to $BG$: an
approximation or the simplicial construction. 

The first approach is to
approximate $BG$ by algebraic varieties $B_nG$. For a detailed
exposition see \cite{T}. If $G\subset
GL(k)$, then for example we take as $E_nG$ the Stiefel variety of
$k$ independent vectors in $\Ct^{n+k}$. The quotient $B_nG=E_nG/G$
exists and it is an algebraic variety. It can be presented as the
quotient of groups 
$GL(n+k)/H$, where $H= \left(\matrix G& *\\ 0 &
GL(n)\\\endmatrix\right)$. Another candidate for an approximation is
the quotient $B'_nG=GL(n+k)/H'$, where $H'= \left(\matrix G& 0\\ 0 &
GL(n)\cr\endmatrix\right)\simeq G\times GL(n)$. In this case
$E'_nG=GL(n+k)/H'_0$, where $H_0'= \left(\matrix I & 0\\ 0 &
GL(n)\endmatrix\right)$. 

Now we construct the associated bundle. Using the
approximation $B'_nG$ we have $E'_nG\times_GX=GL(n+k)\times_{H'}X$ (where
$GL(n)$ acts trivially on $X$). But even here, with such a simple
construction the resulting space does not have to be a variety.
The spaces $E'_nG\times_G X$ (approximations of the Borel
construction) are
 quotients of $E'_nG\times X$ by $G$ and it can be
shown, that they are algebraic spaces.
 Let us quote \cite{B-B}:
\s
\proclaim{ Theorem 2 } Let $G$ be an algebraic group. Let $E$ and
$X$ be algebraic spaces and let $E$ be a principal locally
isotrivial $G$-fibration. Assume that $E$ is normal. Then
$E\times_GX$ exists in the category of algebraic spaces. If,
moreover $E$ is an algebraic variety, $X$ normal and can be
covered by $G$-invariant open quasi-projective subsets, then
$E\times_GX$ is an algebraic variety. \endproclaim

Note that, by \cite{B-B, Theorem 3}, if $H$ is a connected affine
algebraic group, $G\subset H$ and $X$ normal, then $H\times_GX$ is
an algebraic variety if and only if the assumption about covers is
fulfilled. For the proof of formality one would like to apply the
decomposition theorem (or other methods based on Weil conjectures).
It is clear that either we have to make some assumptions about $X$,
or to generalize the results of \cite{BBD} to the case algebraic spaces.

\s

In order to avoid algebraic spaces it is possible to apply a
simplicial construction of $BG$, see \cite{D2, 6.1}. Then $(BG)_\x$ is a
simplicial variety: $(BG)_n=G^n$. The associated Borel construction is
also a simplicial variety: $(EG\times_GX)_n=G^n\times X$. The
topological realization of the Borel construction coincides  with 
the
realization of the quotient stack $X/G$ as defined in \cite{Si, \S8}.
\s

\remark{ Remark 3} Applying any of the described construction we obtain a 
Hodge structure on
$H^*(BG)$ which is pure, see \cite{D2, 9.1}. If we perform the same
constructions for groups over a finite field, we also obtain a pure
structure with respect to the action of the Frobenius authomorphism.
\endremark

\remark{Remark 4} Suppose the base field is
$\Ft_q$. Then the action of the Frobenius authomorphism $Fr^q:x\mapsto
x^q$ induces an operation on  $H^*_{et}(BG;\Qt_\ell)$,  which  is 
just the Adams
operation $\psi^q$; compare \cite{AM; Sul, Ch.5 \S4}.\endremark

\head 3. Degeneration of the spectral sequence without the assumption
of projectivity \endhead

Suppose $X$ is a complete algebraic variety. We replace $X$ by
its normalization. It does not change the intersection cohomology.
Assume that $X$ can be covered by $G$-invariant open quasi-projective
subsets, then the space $E_nG\times_GX$ (or $E'_nG\times_GX$) is an
algebraic variety by Theorem 2. The projection
$$E_nG\times_GX\rightarrow B_nG$$ is a proper algebraic morphism. \s

\proclaim{ Corollary 5} If the normalization of $X$ can be covered
by $G$-invariant quasiprojective open subsets, then the spectral
sequence (*) converging to $IH^*(E_nG\times_GX)$ degenerates at the $E_2$
term.\endproclaim

\noindent{\it Proof.} The Corollary follows from the decomposition
theorem of \cite{BBD}.\qed
\s
\remark{ Remark 6}
We would be able to remove the assumption about
covers, provided that we proved the decomposition theorem for
separable algebraic spaces. 
The statement of \cite{BBD, 6.2.5} seems to remain valid. To prove
it one would have to develop the
formalism of mixed sheaves for
algebraic spaces over finite fields.\endremark

\head 4. The end of the proof \endhead

If $X$
can be covered by $G$-invariant quasiprojective open subsets, then
$EG\times_GX$ can be approximated by algebraic varieties. The
formality follows from Corollary 5.
For arbitrary complete $X$ there is a way of going around the problem
with $EG\times_GX$ not being a variety. First we consider the maximal
torus $T\subset G$. The classifying space for the torus can be
realized as a product of infinite projective spaces.  Then the
approximation of $ET\rightarrow BT$ is locally Zariski trivial. Thus
$ET\times_TX$ is a limit of algebraic varieties. We can apply directly
the techniques of \cite{BBD} to deduce that the spectral sequence (*)
degenerates for $T$.
The restriction to the
maximal torus induces an 
injection $H^*(BG)\hookrightarrow H^*(BT)$ and an injection of the
$E_2$ terms of the spectral sequences (*). Therefore, since the
differentials vanish in the $E_2$ term for $T$, they vanish
in the $E_2$ term for $G$. Then $E_2=E_3$ and we continue our
reasoning. We find that all the differentials in the spectral
sequence for $G$ vanish. The reduction to the maximal torus was
described in \cite{Gi, \S4, prof of 3.3}.
\s
\remark{Remark 7}
There is another solution of the problem with Borel
construction. We apply the simplicial model of $EG$.
We reduce the base ring to a finite field. Now we follow Deligne's suggestion:
the Galois group acts on the spectral sequence (*). Since
$H^*(BG)$ and $IH^*(X)$ are pure,
$E^{pq}_2$ is pure of weight $p+q$. All the differentials in the
spectral sequence vanish because they mix weights. Note that this
proof works for noncomplete algebraic varieties, provided that
$IH^*(X)$ is pure. One can generalize our proof further and consider
cohomology with coefficients in an arbitrary `sheaf' $\Cal F$, an
object of the derived
category. In this
case a suitable extension of $\Cal F$ on the simplicial variety
$(EG\times_GX)_\x$ 
has to be chosen to make sense of the equivariant cohomology. 
The extension $\Cal F_G$ should be an equivariant sheaf in the sense
of \cite{BL} 
and it should be defined over a finitely generated subring $R\subset \Ct$.
Moreover $\Cal F_G$ should have enough good reductions to finite
fields, so that $H_{G_\Ct}^*(X_R\otimes \Ct;\Cal F_{G_\Ct})\simeq 
H_{G_{\overline\Ft_q}}^*(X_R\otimes \overline\Ft_q;\Cal
F_{G_{\overline\Ft_q}})$ for some 
residue field $\Ft_q=R/\frak m$.
The
conclusion is the same: purity of $H^*(X;\Cal F)$ implies formality.
\endremark

\head 5. Algebra structure \endhead

Theorem 1 shows that the $H^*(BG)$-module structure of $H_G^*(X)$
caries no information about the action. The algebra structure of
$H^*_G(X)$ is more interesting, but hard to recover. It involves some
characteristic classes of the action.

\s\remark{Example 8} Let $X=\Pt^N$ and let $G=GL(N+1)$ with the obvious
action on $X$ (as in the Step 2 of \S1). The effect of the Borel
construction is the projectivization of the tautological bundle over
$BG$. The cohomology algebra is well known:
$$H^*_G(X)=\Qt[c_1,c_2,\dots,c_{N+1},h]/
(c_{N+1}-c_Nh+\dots+c_1(-h)^N+(-h)^{N+1}).$$
\endremark

\head 6. Dual structure \endhead

Due to Koszul duality \cite{GKM; MW} the nonequivariant cohomology of $X$ is 
equipped
with the structure of a module over $H_*(G)$: for $\lambda\in H_k(G)$
and $x\in H^l(X)$ we set $\lambda\cdot x=\lambda\smallsetminus\mu^*(x)\in 
H^{l-k}(X)$,
where $\mu:G\times X\rightarrow X$ is the action and $\smallsetminus$ is the 
slant
product $\smallsetminus:H_k(G)\otimes H^l(G\times X)\rightarrow H^{l-k}(X)$.
According to \cite{GKM, 9.3} if the spectral sequence (*) degenerates,
then the cohomology of $X$ is trivial module over $H_*(G)$, even on
the level of the derived category of $H_*(G)$-modules. Therefore
for any $\lambda$ with $\roman{deg}\,\lambda>0$ we have
$\lambda\cdot x=0$, i.e.:
\proclaim{Corollary 9} The structure of $H_{>0}(G)$-module on $IH^*(X)$
is trivial for a complete variety $X$ acted by a connected group $G$.
\endproclaim
 This generalizes a result of Deligne
\cite{D2, 9.1.7}.

\enddocument
\end